%

\magnification=\magstep1
\input amstex
\documentstyle{amsppt}
\pagewidth{6.5truein}
\pageheight{8.9truein}
\ifx\refstyle\undefinedZQA\else\refstyle{C}\fi
\loadbold

\redefine\qed{ \null\nobreak\hfill$\blacksquare$}
\define\restrict{{\restriction}}
\define\nullset{\varnothing}
\define\nullseq{{\langle\rangle}}
\define\Q{{\bold{Q}}}
\define\covmeag{\text{\rm cov}(\bold K)}
\define\SIGMA{{\boldsymbol\Sigma}}
\define\PI{{\boldsymbol\Pi}}
\define\DELTA{{\boldsymbol\Delta}}
\define\Gdelta{\PI^0_2}
\define\Baire{{{}^\omega\omega}}
\define\Seq{{{}^{{<}\omega}\omega}}
\define\Tree{{\Cal T}}
\define\GTree{\Tree'}
\define\PP{P}
\define\setdiff{\backslash}
\define\forces{\Vdash}
\define\rank#1{{\text{\rm rk}_{#1}}}

\define\vDouwen{1}
\define\Harr{2}
\define\Hurewicz{3}
\define\Kuratowski{4}
\define\Larman{5}
\define\MauBP{6}
\define\MauPTA{7}
\define\MauSSTP{8}
\define\Mosch{9}
\define\Solecki{10}
\define\Steel{11}
\define\Stern{12}

\topmatter
\title On disjoint Borel uniformizations\endtitle
\author Howard Becker and Randall Dougherty\endauthor
\affil University of South Carolina\\Ohio State University\endaffil
\date March 16, 1995 \enddate
\thanks The authors were supported by NSF grant number DMS-9158092; the
second author was also supported
by a fellowship from the Sloan Foundation.\endthanks
\address Department of Mathematics, University of South Carolina,
Columbia, SC 29208\endaddress
\email becker\@cs.scarolina.edu \endemail
\address Department of Mathematics, Ohio State University,
Columbus, OH 43210\endaddress
\email rld\@math.ohio-state.edu \endemail
\abstract
Larman showed that any closed subset of the plane with uncountable
vertical cross-sections has $\aleph_1$ disjoint Borel uniformizing sets.
Here we show that Larman's result is best possible: there exist
closed sets with uncountable cross-sections which do not have
more than $\aleph_1$ disjoint Borel uniformizations, even if
the continuum is much larger than $\aleph_1$.  This negatively
answers some questions of Mauldin.  The proof is based on a result
of Stern, stating that certain Borel sets cannot be written as
a small union of low-level Borel sets.
The proof of the latter result uses Steel's method of forcing with
tagged trees; a full presentation of this method, written in
terms of Baire category rather than forcing, is given here.
\endabstract
\endtopmatter

\document

Let $I$ be the unit interval $[0,1]$.  It is well known that there
exist Borel sets $B \subseteq I \times I$ such that all cross-sections
$B_x = \{y\colon (x,y) \in B\}$ are nonempty but there does not
exist a Borel uniformization of~$B$ (a Borel set $U \subseteq B$ such that
for every~$x$ there is a unique~$y$ such that $(x,y) \in U$; this
can also be viewed as a Borel function from~$I$ to~$I$ which selects
a point from each cross-section).  On the other hand, in some cases
(e.g., if the cross-sections~$B_x$ are all $\sigma$\snug-compact
or all non-meager), one can prove that a Borel uniformization exists.
See Moschovakis \cite{\Mosch} for these results.

If all of the cross-sections $B_x$ are uncountable, then it is natural
to ask whether one can find a large number of disjoint Borel uniformizations
of~$B$.  Larman~\cite{\Larman} has shown that, if the sets~$B_x$
are all uncountable and closed (or just $\DELTA^0_2$), then one
can always find $\aleph_1$~disjoint Borel uniformizations of~$B$.
The main purpose of the present paper is to show that the $\aleph_1$
in Larman's result is best possible.

\proclaim{Theorem 1} There is a closed set $B \subseteq I \times I$
such that all cross-sections $B_x = \{y\colon (x,y) \in B\}$ are
uncountable but there do not exist uncountably many disjoint Borel
uniformizations of~$B$ whose ranks (as Borel functions from~$I$ to~$I$)
are bounded below~$\omega_1$.
\endproclaim

Hence, there cannot exist $\aleph_2$ disjoint Borel uniformizations of
this set~$B$ (since one would be able to choose $\aleph_2$ of them with
the same Borel function rank).  So, unless the Continuum Hypothesis is
true, there do not exist continuum many disjoint Borel uniformizations
of~$B$.  This answers a question raised by Mauldin~\cite{\MauSSTP}.

It also follows from Theorem~1 that there do not exist uncountably many
Borel measurable selector functions of bounded Borel rank for the space
$K(I)$ of nonempty compact subsets of~$I$ which select distinct points
within any uncountable compact set (since, as noted in Mauldin
\cite{\MauSSTP}, such selector functions could be applied to the
cross-sections of~$B$ to get disjoint Borel uniformizations of~$B$).
Hence, one cannot find $\aleph_2$ (or $2^{\aleph_0}$ if CH fails) Borel
measurable selector functions for this space which select distinct
points within any uncountable compact set.  (Mauldin~\cite{\MauSSTP}
had shown that one can find $\aleph_1$ such functions.) This settles
Problem~5.1 from Mauldin~\cite{\MauPTA}.

The main step in the proof of Theorem~1 is the following result.
Let $\covmeag$ be the least cardinal~$\kappa$ such that a
perfect Polish space can be expressed as a union of $\kappa$
meager sets.  (It does not matter which perfect Polish space
is used to define $\covmeag$, because any such space has a
comeager subset homeomorphic to the Baire space.)
Clearly $\aleph_1 \le \covmeag \le 2^{\aleph_0}$.

\proclaim{Theorem 2 (Stern)} For any $\alpha < \omega_1$, there is a Borel
subset of the Baire space~$\Baire$ which cannot be expressed as the union of
fewer than $\covmeag$ $\PI^0_\alpha$ sets. \endproclaim

The proof of this, described in Stern~\cite{\Stern} (although the
result is not stated as generally there), uses Steel's method of
forcing with tagged trees.  Actually, Stern combines this method with an
analysis of the Borel ranks of collections of
well-founded trees to produce a stronger result: for any
$\alpha < \omega_1$, any Borel set
which is a union of fewer than $\covmeag$ $\SIGMA^0_\alpha$ sets
must itself be $\SIGMA^0_\alpha$.  The weaker version above (which
was rediscovered independently by the authors) suffices
for the application here.  

The method of Steel forcing is
presented in Harrington~\cite{\Harr} and Steel~\cite{\Steel}; we will
give another presentation here, in terms of Baire
category rather than forcing.

Solecki~\cite{\Solecki} has recently given a different proof of
Stern's results, using effective descriptive set theory.

\proclaim{Corollary 3} A complete analytic or coanalytic set in an
uncountable Polish space cannot be written as a union of fewer than
$\covmeag$ Borel sets with ranks bounded below~$\omega_1$. \endproclaim

(It is well known that any analytic or coanalytic set can be
written as a union of $\aleph_1$ Borel sets \cite{\Mosch}.)

\demo{Proof} Let $X$ be an analytic (or coanalytic) subset of~$\Baire$
which is complete for analytic (coanalytic) subsets of~$\Baire$ using
continuous maps.  It will suffice to show that
$X$~cannot be written as a union of fewer than
$\covmeag$ Borel sets with ranks bounded below~$\omega_1$, because
if $Y$ were a complete analytic (coanalytic) set which could
be written as such a union, then one could fix a Borel map
reducing $X$ to~$Y$ and take preimages of the Borel sets of bounded rank
with union~$Y$ to get Borel sets of bounded rank with union~$X$.
Now, for any $\alpha < \omega_1$, we can find a Borel set $W \subseteq
\Baire$ as in Theorem~2.  Let $g\colon\Baire\to\Baire$ be a continuous
map reducing $W$ to~$X$.  Then $X$ cannot be a union of
fewer than $\covmeag$ $\PI^0_\alpha$ sets, because, if it were,
one could take preimages under~$g$ to get fewer than
$\covmeag$ $\PI^0_\alpha$ sets
with union~$W$, which is impossible.  Since $\alpha$ was arbitrary,
we are done. \qed\enddemo

Actually, one can get a slightly stronger result.

\proclaim{Corollary 4} In any uncountable Polish space, there exist two
disjoint coanalytic sets which cannot be separated by a set which is
a union of fewer than $\covmeag$ Borel sets of bounded rank. \endproclaim

\demo{Proof} Since all uncountable Polish spaces are Borel isomorphic,
it will suffice to work in the space~$(\Baire)^3$.  We follow the usual
construction of a universal pair of disjoint coanalytic sets:  Let
$U$ be a universal coanalytic set in~$(\Baire)^2$ (i.e., all
coanalytic subsets of $\Baire$ occur as cross-sections~$U_x$),
let $C = \{(x,y,z)\colon (x,z) \in U\}$
and $D = \{(x,y,z)\colon (y,z) \in U\}$,
and apply the reduction principle for coanalytic sets to get
disjoint coanalytic sets $C' \subseteq C$ and $D' \subseteq D$
such that $C' \cup D' = C \cup D$.  Now, for any $\alpha < \omega_1$,
let $B$ be the Borel set obtained from Theorem~2, and find $x$ and~$y$
such that $U_x = B$ and $U_y = \Baire\setdiff B$.  Then
$C_{x,y} = B$ and $D_{x,y} = \Baire\setdiff B$, so
$C'_{x,y} = B$ and $D'_{x,y} = \Baire\setdiff B$.
Hence, $C'_{x,y}$ and~$D'_{x,y}$ cannot be separated by a union
of fewer than $\covmeag$ $\PI^0_\alpha$ sets, so
$C'$ and~$D'$ cannot either. \qed\enddemo

Of course, the preceding results say little if $\covmeag = \aleph_1$
(e.g., if CH holds).  However, under Martin's Axiom, the union of fewer
than $2^{\aleph_0}$ meager sets is meager, so $\covmeag = 2^{\aleph_0}$
and these results are more interesting.

In order to prove Theorem~1, we will need to use Corollary~4 to rule
out separating sets that are the union of $\aleph_1$ Borel sets
of bounded rank.  This can be done directly if $\covmeag > \aleph_1$;
if $\covmeag = \aleph_1$, then we will need to do a forcing and
absoluteness argument.

The closed set for Theorem~1 will be obtained from a construction
given in Mauldin \cite{\MauBP, Example~3.2}.  The construction
uses the following well-known fact, proved by methods probably
due to Hurewicz~\cite{\Hurewicz}.

\proclaim{Lemma 5 \rm (Hurewicz?)} For any analytic set $A \subseteq I$,
there is a closed set $B \subseteq I \times I$ such that: if $x \in A$,
then $B_x$ is uncountable; if $x \notin A$, then $B_x \subseteq \Q$.
\endproclaim

\demo{Proof} Since $A$ is analytic, there is a closed set $C \subseteq
I \times \Baire$ whose projection to~$I$ is~$A$.  Define
$C' \subseteq I \times \Baire \times \Baire$ so that
$(x,y,z) \in C'$ iff $(x,y) \in C$.  Then, if $x \in A$, there
are uncountably many $(y,z) \in \Baire \times \Baire$ such that
$(x,y,z) \in C'$; if $x \notin A$, then there is no such $(y,z)$.
But $\Baire \times \Baire$ is homeomorphic to~$\Baire$, which is homeomorphic
to the set of irrationals in~$I$; let $f$ be a homeomorphism from
$\Baire \times \Baire$ to the irrationals in~$I$.  Let $B$ be the
closure in $I \times I$ of the set $\{(x,f(y,z))\colon (x,y,z) \in C'\}$;
then $B$ has the desired properties.  \qed\enddemo

\demo{Proof of Theorem 1} Let $D_1$ and~$D_2$ be disjoint inseparable
coanalytic subsets of~$I$ as given by Corollary~4, and let $A_1 =
I\setdiff D_1$ and $A_2 = I \setdiff D_2$; then $A_1 \cup A_2 = I$.
Let $B_1$ and~$B_2$ be closed sets in $I \times I$ obtained by applying
Lemma~5 to $A_1$ and~$A_2$.  Apply linear mappings to the second
coordinate to compress $B_1$ and~$B_2$ to sets $\hat B_1 \subseteq I
\times [0,1/3]$ and $\hat B_2 \subseteq I \times [2/3,1]$ with the same
properties.  Now let $B = \hat B_1 \cup \hat B_2$.  Since $A_1 \cup A_2
= I$, all cross-sections~$B_x$ are uncountable.  It remains to show
that $B$ does not have $\aleph_1$ disjoint Borel uniformizations of
bounded rank.

First, let us assume that $\covmeag > \aleph_1$.  Suppose that we have
a collection $\{u_\gamma \colon \gamma < \omega_1\}$ of pairwise
disjoint functions from~$I$ to~$I$ each uniformizing~$B$, whose Borel
function ranks are bounded by some fixed $\alpha < \omega_1$.  For each
$\gamma$, let $E_\gamma$ be the set of $x \in I$ such that
$u_\gamma(x)$ is an irrational number greater than~$1/2$.  Since the
set of irrational numbers above~$1/2$ is~$\Gdelta$, each set~$E_\gamma$
is $\PI^0_{\alpha+1}$.  Now, if $x \in D_2$, then $x \notin A_2$, so
$B_x$ contains no irrationals above~$1/2$, so $x \notin E_\gamma$ for
all $\gamma$; if $x \in D_1$, then $x \notin A_1$, so $B_x$ contains no
irrationals below~$1/2$, and since the values $u_\gamma(x)$ for $\gamma
< \omega_1$ are distinct, only countably many of them can be rational,
so $x \in E_\gamma$ for all but countably many~$\gamma$.  Therefore,
the set $\bigcup_{\gamma < \omega_1} E_\gamma$ is a union of
$\aleph_1$~Borel sets of bounded rank which separates $D_1$ from~$D_2$;
since $D_1$ and~$D_2$ were obtained from Corollary~4, and $\aleph_1 <
\covmeag$, we have a contradiction.

Now let us drop the assumption that $\covmeag > \aleph_1$.  Suppose
that we have disjoint functions $u_\gamma$ for $\gamma < \omega_1$
as above.  Fix Borel codes of the appropriate ranks for the functions
$u_\gamma$ and the set~$B$.  By going through the details of the
construction of~$B$, one can check that one obtains the same Borel
code for~$B$ no matter what transitive model of set theory one is working
in.  Now, using the current universe as the ground model, construct
a generic extension with the same~$\omega_1$ in which Martin's Axiom
plus $\neg \text{CH}$ holds.  (Collapse some cardinals above $\aleph_1$
in order to make $2^{\aleph_0} = \aleph_1$ and $2^{\aleph_1} = \aleph_2$,
and then do the standard c.c.c.~forcing iteration to get
$\text{MA} + 2^{\aleph_0}{=}\aleph_2$.)  All of the properties
we assumed about the functions~$u_\gamma$, including the
property of being a function with domain~$I$, are easily seen to be
$\PI^1_2$ assertions about the Borel codes (which are only used
one or two at a time), so, by the Shoenfield absoluteness theorem,
these codes define functions in the generic extension which satisfy
the same assertions.  But $\covmeag > \aleph_1$ holds in the extension,
so we get a contradiction as in the preceding paragraph. \qed\enddemo

It now remains to give the proof of Theorem 2. 

\demo{Proof of Theorem 2}
Let~$\Tree$ be the space of trees on~$\omega$, viewed as a
(closed) subspace of the space of subsets of~$\Seq$ with the usual
Cantor topology  (which in turn is homeomorphic to a closed subspace of
the Baire space).  For any tree $T \in \Tree$, define the rank function
$\rank T\colon T \to \omega_1\cup\{\infty\}$ by the following
condition: for any $s \in T$, $\rank T(s)$ is an ordinal iff $\rank
T(sn)$ is an ordinal for all immediate successors~$sn$ of~$s$ which are
in~$T$, and in this case $\rank T(s)$ is the least ordinal greater than
all of the ordinals $\rank T(sn)$.  So $\rank T(s) = 0$ iff $s$ is a
leaf of~$T$, and $\rank T(s) = \infty$ iff $T$ is not well-founded
below~$s$.

We will work with a slightly restricted set of trees:
let $\GTree$ be the set of $T \in \Tree$ such that the null
sequence~$\nullseq$ is in~$T$ and, for any sequence~$s$, either all
immediate successors~$sn$ ($n \in \omega$) of~$s$ are in~$T$ or none of
them are.  Clearly $\GTree$ is closed in~$\Tree$.

For any $\beta < \omega_1$, let $R_\beta$ be the set of trees $T \in \GTree$
such that, for any $s \in T$, $\rank T(s)$ is either~$\infty$ or less
than $\beta$.  Equivalently, $R_\beta$ is the set of trees which have
no nodes of rank exactly~$\beta$.  Since an easy induction on $\gamma$
shows that $\{T \in \Tree\colon \rank T(s) = \gamma\}$ is Borel for
any~$s$ and any $\gamma < \omega_1$, the sets~$R_\beta$ are all
Borel.  We will prove Theorem~2 by showing that, if $\beta \ge
\omega\cdot\alpha$, then $R_\beta$ is not a union of fewer than
$\covmeag$ $\PI^0_\alpha$ sets.

Define a {\it tagged tree} to be a pair $(T,H)$ where $T \in \Tree$
and $H$ is a function from $T$ to $\omega_1 \cup \infty$ such that,
for any $s \in T$ and $s' \subset s$, we have $H(s') > H(s)$ (where $\infty$
is defined to be greater than any ordinal and greater than itself).
For example, if $T \in \Tree$, then $(T,\rank T)$ is a tagged tree,
and so is $(T',\rank T \restrict T')$ for any subtree~$T'$ of~$T$.
We write $(T,H) \subseteq (T',H')$ when $T \subseteq T'$ and $H
\subseteq H'$.  For $\beta < \omega_1$, a {\it $\beta$\snug-tagged tree}
is a tagged tree $(T,H)$ such that $H \colon T \to \beta\cup\{\infty\}$.
A $\beta$\snug-tagged tree can be viewed as a subset of
$\Seq \times (\beta \cup \{\infty\})$; the set~$\Tree_\beta$ of all
$\beta$\snug-tagged trees is a $\Gdelta$ subset of the space of
subsets of $\Seq \times (\beta \cup \{\infty\})$, so, with the
inherited topology, it is itself a Polish space by Alexandrov's
Theorem~\cite{\Kuratowski}.

Let $\PP_\beta \subseteq \Tree_\beta$ be the set of finite
$\beta$\snug-tagged trees.  Now define a new topology on~$R_\beta$,
to be called the $\beta$\snug-topology, with basis consisting of the
sets $$N^\beta_p = \{T \in R_\beta\colon p \subseteq (T,\rank T)\}$$
for $p \in \PP_\beta$.  Then the $\beta$\snug-topology is a
Polish topology on~$R_\beta$.  To see this, let $S$ be the set
of $(T,H) \in \Tree_\beta$ such that $T \in \GTree$ and
$H$ satisfies the recursive
definition of~$\rank T$; then it is easy to check that $S$ is
$\Gdelta$ in~$\Tree_\beta$ and hence Polish by Alexandrov's Theorem.
It is not hard to show that the projection $(T,H) \mapsto T$ is
a homeomorphism from~$S$ to $R_\beta$ with the $\beta$\snug-topology.
(One needs the fact that, for any $s \in \Seq$, the set $\{T \in R_\beta
\colon s \notin T\}$ is open in the $\beta$\snug-topology; this set
can in fact be written as the union of~$N^\beta_p$ for those
$p = (t,h)$ such that $h(s') = 0$ for some $s' \subset s$, because
we have restricted ourselves to trees in~$\GTree$.)
So the $\beta$\snug-topology is Polish, and includes the original
topology on~$R_\beta$ as a subspace of~$\Tree$.

For any set~$A\subseteq\Tree$ and any $p \in \PP_\beta$, define
$p\forces_\beta A$ to mean that $A \cap N^\beta_p$ is comeager
in~$N^\beta_p$ under the $\beta$\snug-topology.  Easily, if $p
\subseteq q$, then $p\forces_\beta A$ implies $q\forces_\beta A$; if $A
\subseteq B$, then $p\forces_\beta A$ implies $p\forces_\beta B$; and
$p \forces_\beta \bigcap_{n=0}^\infty A_n$ if and only if $p
\forces_\beta A_n$ for all~$n$.  Furthermore, if $A \cap R_\beta$ has
the Baire property in the $\beta$\snug-topology, then
$p \not\forces_\beta A$ if and only if there is $q \supseteq p$
such that $q \forces_\beta {-A}$.  In particular, this is true
whenever $A$ is a Borel subset of $\Tree$, since then $A \cap R_\beta$
is Borel in~$R_\beta$ under the inherited topology and hence
under the $\beta$\snug-topology as well.

For example, let $p_0$ be the tagged tree $(\{\nullseq\},h)$
where $h(\nullseq) = \infty$; then $R_\beta \setdiff R_{\beta'}$ is
$\beta$\snug-open dense in~$N^\beta_{p_0}$ for any $\beta' < \beta$
(since any $p \supseteq p_0$ in $\PP_\beta$ can be extended by adding
a new sequence of length~$1$ to the tree with tag~$\beta'$),
so $p_0 \forces_\beta R_\beta \setminus \bigcup_{\beta' < \beta}
R_{\beta'}$.

If $(t,h)$ and $(t',h')$ are finite tagged trees and $\alpha$ is
an ordinal, define $(t,h) \sim_\alpha (t',h')$ to mean that
$t = t'$ and, for any $s \in t$, if either of $h(s)$ and~$h'(s)$ is
an ordinal less than~$\alpha$, then $h(s) = h'(s)$.

The following lemma is known as the Retagging Lemma.

\proclaim{Lemma 6 \rm (Steel)} If $\alpha \ge 1$ is a countable ordinal,
$\beta_1,\beta_2 \ge \omega\cdot\alpha$, $p_1 \in \PP_{\beta_1}$,
$p_2 \in \PP_{\beta_2}$, and $p_1 \sim_{\omega\cdot\alpha} p_2$, then, for
any $\PI^0_\alpha$ set $A \subseteq \Tree$, $p_1 \forces_{\beta_1} A$ if
and only if $p_2 \forces_{\beta_2} A$. \endproclaim

\demo{Proof}
Say $p_1 = (t,h_1)$ and $p_2 = (t,h_2)$.
We will proceed by induction on $\alpha$.

For $\alpha = 1$, suppose that $p_1 \not\forces_{\beta_1} A$; we will
show that $p_2 \not\forces_{\beta_2} A$.  (Of course, the reverse
implication is identical.)  Let $T$ be a tree in~$N^{\beta_1}_{p_1}$
which is not in~$A$.  Since $A$ is~$\PI^0_1$, $-A$ is open, so there
exist finitely many sequences $s_1,\dots,s_m \in T$ and
$s'_1,\dots,s'_k \notin T$ such that any tree containing all of the
sequences $s_i$ and none of the sequences $s'_i$ is in~$-A$.  For each
$i \le k$, let $s''_i$ be the longest initial segment of $s'_i$ that is
in~$T$; then, since $T \in \GTree$, each~$s''_i$ is a leaf of~$T$
(i.e., $\rank T(s''_i) = 0$).  Now let $\tau \supseteq t$ be a finite
subtree of~$T$ containing all of the sequences $s_i$ and~$s''_i$, and
let $q = (\tau, h)$ where $h = \rank T\restrict \tau$; then $q
\supseteq p_1$ and we have not only $q \forces_{\beta_1} {-A}$, but also
$\bar q \forces_{\bar\beta} {-A}$ whenever $\bar q \in \PP_{\bar\beta}$
and $\bar q \sim_1 q$ (since, in any tree in~$N^{\bar \beta}_{\bar
q}$, all sequences $s_i$ would be nodes and all sequences $s''_i$
would be leaves).  Let $M$ be a natural number greater than all natural
numbers occurring as tags in~$q$ or~$p_2$, and let $\gamma_0 < \dots <
\gamma_{n-1}$ list the infinite ordinals occurring as tags in~$q$.
Let $L$ be the largest of the lengths of the sequences in~$\tau$.  Now
define $\hat h \colon \tau \to \beta_2 \cup \{\infty\}$ as follows:
$$\hat h(s) = \cases h_2(s) &\text{if $s \in t$},\\
h(s) &\text{if $s \notin t$ and $h(s) < \omega$},\\
M+j &\text{if $s \notin t$ and $h(s) = \gamma_j$},\\
M+n+L-\text{len}(s) &\text{if $s \notin t$ and $h(s) = \infty$}. \endcases$$
Then, since $p_1 \sim_\omega p_2$, it is easy to check that $(\tau,\hat
h)$ is a valid $\beta_2$\snug-tagged tree extending $p_2$ and
$(\tau,\hat h) \sim_1 q$.  Hence, $(\tau,\hat h) \forces_{\beta_2}
{-A}$, so $p_2 \not\forces_{\beta_2} A$, as desired.

Now suppose $\alpha > 1$, and write~$A$ as a countable intersection of
sets~$A_k$, each of which is~$\SIGMA^0_{\alpha'}$ for some $\alpha' <
\alpha$ (which may vary with~$k$).  Suppose that $p_1
\not\forces_{\beta_1} A$; we will show that $p_2 \not\forces_{\beta_2}
A$.  There must be a~$k$ such that $p_1 \not\forces_{\beta_1} A_k$,
and hence $q \forces_{\beta_1} {-A_k}$ for some $q = (\tau,h)\supseteq
p_1$.  Fix $\alpha' < \alpha$ such that $A_k$ is~$\SIGMA^0_{\alpha'}$ and
hence $-A_k$ is~$\PI^0_{\alpha'}$.  Arguing as before, let $M$ be a
natural number such that $\omega\cdot\alpha'+M$ is greater than all
ordinals below $\omega\cdot \alpha' +\omega$ occurring as tags in~$q$
or~$p_2$, and let $\gamma_0 < \dots < \gamma_{n-1}$ list the ordinals at or
above $\omega\cdot\alpha'+\omega$ occurring as tags in~$q$.  Let $L$ be
the largest of the lengths of the sequences in~$\tau$.  Now define
$\hat h \colon \tau \to \beta_2 \cup \{\infty\}$ as follows:
$$\hat h(s) = \cases h_2(s) &\text{if $s \in t$},\\
h(s) &\text{if $s \notin t$ and $h(s) < \omega\cdot\alpha'+\omega$},\\
\omega\cdot\alpha'+M+j &\text{if $s \notin t$ and $h(s) = \gamma_j$},\\
\omega\cdot\alpha'+M+n+L-\text{len}(s) &\text{if $s \notin t$ and
   $h(s) = \infty$}. \endcases$$
Then, since $p_1 \sim_{\omega\cdot\alpha} p_2$ and $\omega\cdot\alpha
\ge \omega\cdot\alpha'+\omega$, it is easy to check that $(\tau,\hat
h)$ is a valid $\beta_2$\snug-tagged tree extending $p_2$ and
$(\tau,\hat h) \sim_{\omega\cdot\alpha'} q$.  Hence, by the
induction hypothesis, $(\tau,\hat h) \forces_{\beta_2}
{-A_k}$, so $p_2 \not\forces_{\beta_2} A$, as desired.
\qed\enddemo

We are now ready to show that, if $\beta \ge \omega\cdot\alpha$, then
$R_\beta$ cannot be expressed as the union of fewer than~$\covmeag$
$\PI^0_\alpha$ subsets of~$\Tree$.  Suppose it can.  Then these subsets
cover $N^\beta_{p_0}$, which can be viewed as a Polish space under the
$\beta$\snug-topology, so, by the definition of~$\covmeag$, at least
one of these $\PI^0_\alpha$ sets, say~$W$, must be
$\beta$\snug-nonmeager in~$N^\beta_{p_0}$.  We now have $p_0
\not\forces_\beta {-W}$, so there exists a $q \supseteq p_0$ in
$\PP_\beta$ such that $q \forces_\beta W$.  By the Retagging Lemma, we
have $q \forces_\gamma W$ for any $\gamma > \beta$.  As noted before,
$p_0 \forces_\gamma R_\gamma\setdiff R_\beta$, so $q \forces_\gamma W
\cap (R_\gamma\setdiff R_\beta)$.  But $W \subseteq R_\beta$, so $q
\forces_\gamma \nullset$, which is impossible.  This completes the
proof.
\qed\enddemo

A question raised by these results is: exactly what is the least
cardinal~$\lambda$ such that any analytic set is the union of $\lambda$~Borel
sets of bounded rank?  (The same number of Borel sets of bounded rank
would also suffice to give any coanalytic set or even any $\SIGMA^1_2$
set, since a $\SIGMA^1_2$~set is a union of $\aleph_1$~analytic sets.)
Corollary~3 gives a lower bound of~$\covmeag$ for~$\lambda$, while
Theorem~8.10(e) of van~Douwen~\cite{\vDouwen} implies that
the dominating number~$\frak d$ is an upper bound for~$\lambda$,
since it states that any analytic set is a union of $\frak d$ {\sl
compact} sets.  Solecki~\cite{\Solecki} gives related results.


\Refs

\ref \no \vDouwen \by E. K. van Douwen \paper The integers and topology
\inbook Handbook of set-theoretic topology \eds K. Kunen and
J.~E.~Vaughan \publ North-Holland \publaddr Amsterdam \yr 1984
\pages 111--167 \endref

\ref \no \Harr \by L. Harrington \paper Analytic determinacy and $0^\sharp$
\jour J. Symbolic Logic \vol 43 \yr 1978 \pages 685--693 \endref

\ref \no \Hurewicz \by W. Hurewicz \paper Zur Theorie der analytischen
Mengen \jour Fund. Math. \vol 15 \yr 1930 \pages 4--17 \endref

\ref \no \Kuratowski \by K. Kuratowski \book Topology \vol 1
\publ Academic Press \publaddr New York \yr 1966 \endref

\ref \no \Larman \by D. G. Larman \paper Projecting and uniformizing
Borel sets with $K_\sigma$~sets, II \jour Mathematika
\vol 20 \yr 1973 \pages 233--246 \endref

\ref \no \MauBP \by R. D. Mauldin \paper Borel parametrizations
\jour Trans. Amer. Math. Soc. \vol 250 \yr 1979 \pages 223--234 \endref

\ref \no \MauPTA \bysame \paper Problems in topology arising
from analysis \inbook Open problems in topology \eds J. van Mill and
G.~M.~Reed \publ North-Holland \publaddr Amsterdam \yr 1990
\pages 617--629 \endref

\ref \no \MauSSTP \bysame \paper Some selection theorems and problems
\inbook Measure theory, Oberwolfach 1979 \ed D. K\"olzow \publ
Springer-Verlag \publaddr Berlin \bookinfo Lecture notes in mathematics,
vol. 794 \pages 160--165 \endref

\ref \no \Mosch \by Y. N. Moschovakis \book Descriptive set theory
\publ North-Holland \publaddr Amsterdam \yr 1980 \endref

\ref \no \Solecki \by S. Solecki \paper Decomposing Borel sets and
functions and the structure of Baire class 1 functions \paperinfo
in preparation \endref

\ref \no \Steel \by J. R. Steel \paper Forcing with tagged trees \jour
Ann. Math. Logic \vol 15 \yr 1978 \pages 55--74 \endref

\ref \no \Stern \by J. Stern \paper \'Evaluation du rang de Borel
de certains ensembles \jour C. R. Acad. Sci. Paris S\'er. A--B
\vol 286 \yr 1978 \pages A855--A857 \endref

\endRefs
\enddocument